\documentstyle[11pt]{article}
\def\Bbb R{{\rm \bf R}}
\def\proclaim#1{\vskip2mm{\bf #1}\em}
\def\endproclaim{\em \vskip2mm}
\def\tag#1{\eqno(#1)}
\def\gathered{\begin{array}{c}}
\def\endgathered{\end{array}}
\def\text{\mbox}

\begin{document}

\title {On reconstruction from a partial knowledge of the Neumann-to-Dirichlet operator}
\author{Masaru IKEHATA\footnote{
Department of Mathematics,
Faculty of Engineering,
Gunma University, Kiryu 376-8515, JAPAN}}
\date{}
\maketitle

\begin{abstract}
We give formulae that yield an information about the location of an unknown
polygonal inclusion having unknown constant conductivity inside a known
conductive material having known constant conductivity from a partial knowledge of the Neumann
-to-Dirichlet operator.

%{\bf PartialKnowledgeND.tex}

\end{abstract}

%\tableofcontents

\section{Statement of the result}

We give another application of the method developed in \cite{Ie} and \cite{Ie2} to a special,
however, important version of the inverse boundary value problem formulated by Calder\'on \cite{C}.
It is a mathematical formulation of electrical impedance tomography.

Let $\Omega$ be a two-dimensional bounded domain with smooth boundary.
We consider $\Omega$ an isotropic, electrically conductive material.  Let $\gamma$ be the conductivity of $\Omega$,
and $D$ an open set of $\Omega$ such that $\overline D\subset\Omega$ and $\partial D$ is Lipschitz.
Assume that $\gamma$ takes a positive constant $k_j$ on each connected component $D_j$ of $D$ with $k_j$,
$j=1,\cdots,m\le\infty$ and is equal to $1$ on $\Omega\setminus D$.
In this paper, we always assume that $m<\infty$ and that $\overline{D_j}\cap\overline{D_{j'}}=\emptyset$
if $j\not=j'$.  We call $D$ an inclusion and $\gamma$ the corresponding conductivity.

Let $\nu$ denote the unit outward normal vector field to $\partial(\Omega\setminus\overline D)$.
We prescribe the electric current distribution $g\in L^{\infty}(\partial\Omega)$
satisfying $\displaystyle\int_{\partial\Omega}g=0$.
Consider the elliptic problem
$$\displaystyle
\nabla\cdot\gamma\nabla u=0\,\,\text{in}\,\Omega,\,\,\,\,
\frac{\partial u}{\partial\nu}=g\,\,\text{on}\,\partial\Omega.
\tag {1.1}
$$
The problem (1.1) has a $H^1(\Omega)$ solution $u$ and any solution has the form $u+c$ where
$c$ is a constant.  It is well known that $u$ is H\"older continuous on $\overline\Omega$ \cite{L}.

We consider the following problem.

{\it\noindent Inverse Problem}

Let $P$, $Q$ be two arbitrary distinct points on $\partial\Omega$.  We fix $P$ and $Q$.
Then the map
$$\displaystyle
\Lambda_{\gamma}(P,Q): g\longmapsto u(P)-u(Q)
$$
is well defined where $u$ is a solution to the problem (1.1).
Assume that $D$ and $\gamma$ on $D$ are unknown.  The problem is to find a formula
that yields an information about the location of $D$ from $\Lambda_{\gamma}(P,Q)$.
This map is a partial knowledge of the Neumann-to-Dirichlet operator:
$$\displaystyle
\Lambda_{\gamma}:g\longmapsto u\vert_{\partial\Omega}.
$$

$\Lambda_{\gamma}$ uniquely determines $D$ and $\gamma$.  This is a corollary of Isakov's uniqueness theorem \cite{IS}
which covers a more general case.  Under suitable regularity assumption on $\gamma$,
Nachman \cite{N2} established a reconstruction formula of $\gamma$ itself from the full knowledge of $\Lambda_{\gamma}$.
See the survey paper \cite{U} for several other results.

In \cite{Ie2}, using the method developed in \cite{Ie}, we gave formulae that yield the information
about the convex hull of $D$ from $\{u(y)-u(Q)\,\vert\,y\in\partial\Omega\}$
for fixed $Q$ and $g\not\equiv 0$ provided each connected component of $D$ is a polygon
and $D$ satisfies the condition
$$\displaystyle
\text{diam}\,D<\text{dist}\,(D,\partial\Omega).
\tag {1.2}
$$
This result gives a constructive proof of Friedman-Isakov's uniqueness theorem \cite{FI}.
However, the formulae involve an integral of $g(y)$ and of the measured data $u(y)-u(Q)$.

The strong motivation of our study is to seek formulae that avoid any integration
of the measured data on $\partial\Omega$.
In this paper we give two formulae that yield the information about the convex hull of $D$ from
$\Lambda_{\gamma}(P,Q)$ for fixed $P$ and $Q$ provided that each connected component
of $D$ is a polygon and $D$ satisfies (1.2).

Now we describe the result more precisely.
Let $S^1$ denote the set of all unit vectors of $\Bbb R^2$.
Define
$$\displaystyle
h_D(\omega)
=\sup_{x\in D}\,x\cdot\omega,\,\,\,\,\omega\in S^1.
$$
This function is called the support function of $D$ and the convex hull of $D$
can be reconstructed from this function.

We say that $\omega\in S^1$ is regular with respect to $D$ if the set
$$\displaystyle
\{x\in\Bbb R^2\,\vert\,x\cdot\omega=h_D(\omega)\}
\cap\partial D
$$
consists of only one point.

Note that the set of all unit vectors that are not regular with respect to $D$ is a finite set.
Since $h_D(\,\cdot\,)$ is continuous on $S^1$, the convex hull of $D$ can be reconstructed
from the restriction of $h_D(\,\cdot\,)$ to the set of all unit vectors that are regular with respect to $D$.

The following special harmonic functions play the central role:
$$\displaystyle
v_{\omega}=v_{\omega}(x;\tau)
=e^{\tau x\cdot(\omega+i\omega^{\perp})},\,\,\,\,\tau>0
$$
where $\omega, \omega\in S^1$ and satisfy
$$\displaystyle
\omega\cdot\omega^{\perp}=0,\,\,\text{det}\,(\omega\,\,\omega^{\perp})<0.
$$

Define
$$\displaystyle
g=g_{\omega}(y;\tau)
=\frac{\partial v_{\omega}}{\partial\nu}\vert_{\partial\Omega}.
\tag {1.3}
$$
Note that $g\in C^{\infty}(\partial\Omega)$ and $\displaystyle\int_{\partial\Omega}g=0$.

{\bf\noindent Definition 1.1(Indicator function).}
Define
$$\displaystyle
I_{\omega}(\tau,t)
=e^{-\tau t}
\{\Lambda_{\gamma}(P,Q)-\Lambda_1(P,Q)\}g_{\omega}(\,\cdot\,;\tau),\,\,\,\,\tau>0,\,t\in\Bbb R.
$$

Note that one can rewrite
$$\displaystyle
I_{\omega}(\tau,t)
=\{\Lambda_{\gamma}(P,Q)-\Lambda_1(P,Q)\}
\left(\frac{\partial(e^{-\tau t}v_{\omega})}{\partial\nu}\vert_{\partial\Omega}\right).
$$
Since
$$\displaystyle
e^{-\tau t}v_{\omega}
=e^{\tau(x\cdot\omega-t)}e^{i\tau x\cdot\omega^{\perp}},
$$
the function $e^{-\tau t}v_{\omega}$ has a special property as $\tau\longrightarrow\infty$:
$\vert e^{-\tau t}v_{\omega}\vert$ is exponentially growing in the half space $x\cdot\omega>t$ and
exponentially decaying in $x\cdot\omega<t$;
$e^{-\tau t}v_{\omega}$ is oscillating on the line $x\cdot\omega=t$.

The result is the following two formulae.

\proclaim{\noindent Theorem 1.1.}
Assume that each connected component of $D$ is a polygon and $D$ satisfies (1.2).
Let $\omega$ be regular with respect to $D$.  The formulae
$$\begin{array}{c}
\displaystyle
\{t\in\Bbb R\,\vert\,\lim_{\tau\longrightarrow\infty}I_{\omega}(\tau,t)=0\}
=[h_D(\omega),\,\infty [\\
\\
\displaystyle
\lim_{\tau\longrightarrow\infty}\frac{\log\vert I_{\omega}(\tau,t)\vert}{\tau}
=h_D(\omega)-t,\,\,\forall t\in\Bbb R,
\end{array}
$$
are valid.

\endproclaim

In particular, this theorem tells us that
$$\displaystyle
h_D(\omega)\tau\approx\log\vert\{\Lambda_{\gamma}(P,Q)-\Lambda_1(P,Q)\}g_{\omega}(\,\cdot\,;\tau)\vert
$$
for $\tau>>1$.  This will be usefull to calculate the approximate value of $h_D(\omega)$ from $\Lambda_{\gamma}(P,Q)-\Lambda_1(P,Q)$.
We will do the test in the future.

A brief outline of the proof of Theorem 1.1 is as follows.
In Section 2 we construct
a solution $\mbox{$\cal D$}=\mbox{$\cal D$}(P,Q;x)$ to the elliptic problem (see Proposition 2.1)
$$\displaystyle
\nabla\cdot\gamma\nabla\mbox{$\cal D$}=0\,\,\text{in}\,\Omega\,\,\,\,
\frac{\partial\mbox{$\cal D$}}{\partial\nu}=\delta_P-\delta_Q\,\,\text{on}\,\partial\Omega
\tag {1.4}
$$
where for $x\in\partial\Omega$, $\delta_x$ denotes the Dirac measure on $\partial\Omega$ concentrated at $x$.
Using this solution, in Section 3 we establish the representation formula of the indicator function:
$$\displaystyle
I_{\omega}(\tau,t)
=e^{-\tau t}\sum_j(k_j-1)\int_{\partial D_j}\mbox{$\cal D$}(P,Q;\,\cdot\,)\frac{\partial v_{\omega}(\,\cdot\,;\tau)}{\partial\nu}.
$$
From this and the trivial identity
$$\displaystyle
I_{\omega}(\tau,t)
=e^{\tau(h_D(\omega)-t)}I_{\omega}(\tau,h_D(\omega)),
\tag {1.5}
$$
we know that it suffices to study the asymptotic behaviour as $\tau\longrightarrow\infty$ of the oscillatory integral
$$\displaystyle
e^{-\tau t}\sum_j(k_j-1)\int_{\partial D_j}u\frac{\partial v_{\omega}(\,\cdot\,;\tau)}{\partial\nu}
$$
for $u=\mbox{$\cal D$}(P,Q;x)$.  However, we have already encountered this type of integral \cite{Ie, Ie2}.
we see that this integral decays algebraically as $\tau\longrightarrow\infty$ provided (1.2) (see Lemma 3.1).
Then from (1.5) we obtain the two formulae in Theorem 1.1.

\section{Preliminaries}

First we construct a solution to the problem (1.4).
The construction is similar to that of a solution to the crack problem described in Appendix D, D.1 of \cite{AD}.
Set
$$\begin{array}{c}
\displaystyle
V=V(P,Q;x)
=-\frac{1}{\pi}
(\log\vert P-x\vert-\log\vert Q-x\vert),\,\,x\in\Bbb R^2,\\
\\
\displaystyle
\Psi=\Psi(P,Q;x)
=-\frac{1}{\pi}
\left(\frac{y-P}{\vert y-P\vert^2}
-\frac{y-Q}{\vert y-Q\vert^2}\right)\cdot\nu(y),\,\,y\in\partial\Omega\setminus\{P,Q\}.
\end{array}
$$
Since we assumed that $\partial\Omega$ is smooth, $\Psi(P,Q;\,\cdot\,)\in L^{\infty}(\partial\Omega)$.
Moreover, we know that
$$\displaystyle
\int_{\partial\Omega}\Psi(P,Q;y)d\sigma(y)=0.
$$
These are well known facts in the potential theory.
From the assumption $\text{supp}\,(\gamma-1)\subset\overline D$,
we see that
$$\displaystyle
(\gamma-1)\nabla V\in L^{\infty}(\Omega).
$$
Therefore we have the unique $H^1(\Omega)$ solution
$\mbox{$\cal E$}=\mbox{$\cal E$}(P,Q;x)$ to the elliptic problem
$$\begin{array}{c}
\displaystyle
\nabla\cdot\gamma\nabla\mbox{$\cal E$}=-\nabla\cdot(\gamma-1)\nabla V\,\,\text{in}\,\Omega,\\
\\
\displaystyle
\frac{\partial}{\partial\nu}\mbox{$\cal E$}
=-\Psi\,\,\text{on}\,\partial\Omega,\\
\\
\displaystyle
\int_{\partial\Omega}\mbox{$\cal E$}=0.
\end{array}
\tag {2.1}
$$
Define
$$\displaystyle
\mbox{$\cal D$}=\mbox{$\cal D$}(P,Q;x)=V(P,Q;x)+\mbox{$\cal E$}(P,Q;x),\,\,x\in\Omega.
$$
We prove that this satisfies (1.4) in the following sense.

\proclaim{\noindent Proposition 2.1.}
For any $\varphi\in H^1(\Omega)$ that is smooth in a neighbourhood of $\partial\Omega$ the formula
$$\displaystyle
\int_{\Omega}\gamma\nabla\mbox{$\cal D$}(P,Q;x)\cdot\nabla\varphi dx=\varphi(P)-\varphi(Q),
\tag {2.2}
$$
is valid.

\endproclaim

{\it\noindent Proof.}
Note that $\nabla V\cdot\nabla\varphi$ is absolutely integrable in the whole domain because of the regularity 
assumption of $\varphi$.  Therefore, we have
$$\displaystyle
\int_{\Omega}\gamma\nabla\mbox{$\cal D$}\cdot\nabla\varphi dx
=\int_{\Omega}\gamma\nabla V\cdot\nabla\varphi dx+\int_{\Omega}\gamma\nabla\mbox{$\cal E$}\cdot\nabla\varphi dx.
\tag {2.3}
$$
From (2.1) we have
$$\displaystyle
\int_{\Omega}\gamma\nabla\mbox{$\cal E$}\cdot\nabla\varphi dx
=-\int_{\Omega}(\gamma-1)\nabla V\cdot\nabla\varphi dx-\int_{\partial\Omega}\Psi\varphi.
\tag {2.4}
$$
It is well known that
$$\displaystyle
\int_{\Omega}\gamma\nabla V\cdot\nabla\varphi dx
=\int_{\partial\Omega}\Psi\varphi+\varphi(P)-\varphi(Q).
\tag {2.5}
$$
A combination of (2.3)-(2.5) gives (2.2).

\noindent
$\Box$

Using the function $\mbox{$\cal D$}(P,Q;x)$, we obtain a representation formula of $w(P)-w(Q)$ for any
$H^1(\Omega)$ solution $w$ to the problem
$$\displaystyle
\nabla\cdot\gamma\nabla w=-\nabla\cdot F\,\,\text{in}\,\Omega,\,\,\,\,
\frac{\partial}{\partial\nu}w=0\,\,\text{on}\,\partial\Omega
\tag {2.6}
$$
where $F\in L^2(\Omega)$ and satisfy
$$\displaystyle
\text{supp}\,F\subset\overline D.
$$
From this and $\gamma=1$ in a neighbourhood of $\partial\Omega$, $w$ is smooth in a neighbourhood of $\partial\Omega$.

\proclaim{\noindent Proposition 2.2.}  The formula
$$\displaystyle
w(P)-w(Q)
=-\int_DF(x)\cdot\nabla\mbox{$\cal D$}(P,Q;x)dx,
\tag {2.7}
$$
is valid.

\endproclaim

{\it\noindent Proof.}
From (2.2) we have
$$\displaystyle
w(P)-w(Q)=\int_{\Omega}\gamma\nabla w\cdot\nabla\mbox{$\cal D$}dx.
\tag {2.8}
$$
For $\epsilon>0$ define
$$\displaystyle
\mbox{$\cal D$}_{\epsilon}(P,Q;x)
=V_{\epsilon}(P,Q;x)+\mbox{$\cal E$}(P,Q;x)
$$
where
$$\displaystyle
V_{\epsilon}(P,Q;x)
=-\frac{1}{\pi}
\left(\log\sqrt{\vert P-x\vert^2+\epsilon^2}
-\log\sqrt{\vert Q-x\vert^2+\epsilon^2}\right),
\,\,x\in\Bbb R^2.
$$
From the regularity of $w$, Lebesgue's dominated convergence theorem and (2.6) we obtain
$$\displaystyle
\int_{\Omega}\gamma\nabla w\cdot\nabla\mbox{$\cal D$}dx
=\lim_{\epsilon\longrightarrow 0}
\int_{\Omega}\gamma\nabla w\cdot\nabla\mbox{$\cal D$}_{\epsilon}dx
=-\lim_{\epsilon\longrightarrow 0}
\int_{\Omega}F\cdot\nabla\mbox{$\cal D$}_{\epsilon}dx
=-\int_DF\cdot\nabla\mbox{$\cal D$}dx.
$$
This together with (2.8) gives (2.7).

\noindent
$\Box$

\section{Proof of Theorem 1.1}

First we give the representation formula of the indicator function.

\proclaim{\noindent Proposition 3.1.}
The formula
$$\displaystyle
I_{\omega}(\tau,t)
=e^{-\tau t}\sum_j(k_j-1)\int_{\partial D_j}\mbox{$\cal D$}(P,Q;\,\cdot\,)\frac{\partial v_{\omega}(\,\cdot\,;\tau)}{\partial\nu},
\tag {3.1}
$$
is valid.

\endproclaim

{\it\noindent Proof.}
Let $u$ be a $H^1(\Omega)$ solution to the problem (1.1).
Given $g\in L^{\infty}(\partial\Omega)$ satisfying $\displaystyle\int_{\partial\Omega}g=0$
let $v\in H^1(\Omega)$ be a harmonic function satisfying $\displaystyle\frac{\partial v}{\partial\nu}=g$
on $\partial\Omega$.  Then $w=u-v$ satisfies
$$\displaystyle
\nabla\cdot\gamma\nabla w=-\nabla\cdot(\gamma-1)\nabla v\,\,\text{in}\,\Omega,\,\,\,\,
\frac{\partial w}{\partial\nu}=0\,\,\text{on}\,\partial\Omega.
$$
Since $v\in C^{\infty}(\Omega)$ and $\text{supp}\,(\gamma-1)\subset\overline D$, we have
$(\gamma-1)\nabla v\in L^{\infty}(\Omega)$.  Then one obtains from (2.7) that
$$\begin{array}{c}
\displaystyle
(u(P)-u(Q))-(v(P)-v(Q))
=w(P)-w(Q)
=-\int_D(\gamma-1)\nabla v\cdot\nabla\mbox{$\cal D$}(P,Q;x)dx\\
\\
\displaystyle
=-\sum_{j}(k_j-1)\int_{D_j}\nabla v\cdot\nabla\mbox{$\cal D$}(P,Q;x)dx.
\end{array}
\tag {3.2}
$$
Since $\mbox{$\cal D$}(P,Q;x)$ is in $H^1(D_j)$, $v$ is in $H^2(D_j)$ and $\partial D_j$ is Lipschitz, from
Lemma 1.5.3.7 in \cite{Gr} one has
$$\displaystyle
\int_{D_j}\nabla v\cdot\nabla\mbox{$\cal D$}(P,Q;x)dx
=-\int_{\partial D_j}\mbox{$\cal D$}(P,Q;\,\cdot\,)\frac{\partial v}{\partial\nu}.
\tag {3.3}
$$
Note that $v$ is outward to $\partial(\Omega\setminus\overline D)$.
From (3.2), (3.3) and the definition of $\Lambda_{\gamma}(P,Q)$ we obtain
$$\displaystyle
\{\Lambda_{\gamma}(P,Q)-\Lambda_1(P,Q)\}g
=\sum_{j}(k_j-1)\int_{\partial D_j}\mbox{$\cal D$}(P,Q;\,\cdot\,)\frac{\partial v}{\partial\nu}.
\tag {3.4}
$$
Now (3.1) is clear.

\noindent
$\Box$

Equation (3.4) is the representation formula of $\Lambda_{\gamma}(P,Q)-\Lambda_1(P,Q)$.

Now we describe a lemma which is the key for the proof of Theorem 1.1.

\proclaim{\noindent Lemma 3.1.}
Assume that each connected component of $D$ is a polygon and $D$ satisfies (1,2);
$u\in H^1_{\text{loc}}(\Omega)$ satisfies
$$\displaystyle
\int_{\Omega}\gamma\nabla u\cdot\nabla\varphi dx=0
$$
for all $\varphi\in C^{\infty}_0(\Omega)$ and is not a constant function.
Let $\omega$ be regular with respect to $D$ and $v=v_{\omega}(x;\tau)$.  There exist
positive constants $L$ and $\mu$ such that
$$\displaystyle
\lim_{\tau\longrightarrow\infty}
\tau^{\mu}e^{-\tau h_D(\omega)}\left\vert
\sum_j(k_j-1)\int_{\partial D_j}u\frac{\partial v}{\partial\nu}\right\vert=L.
\tag {3.5}
$$

\endproclaim

The proof of this lemma is essentially the same as that of the key lemma in \cite{Ie2}.
We describe only the outline.

{\it\noindent Outline of the proof.}
We have
$$\displaystyle
\int_{\partial D_j}u\frac{\partial v}{\partial\nu}=\int_{\partial D_j}(u-c)\frac{\partial v}{\partial\nu}
\tag {3.6}
$$
for any constant $c$ since $v$ is harmonic in $\Omega$.

From the regularity of $\omega$ we know that the line $x\cdot\omega=h_D(\omega)$ meets $\partial D$
at a vertex $x_0$ of $D$.
Moreover, there exist $j_0$ and $\delta>0$ such that $x_0$ is a vertex of $D_{j_0}$, $h_{D_{j_0}}(\omega)=h_D(\omega)$
and $\cup_{j\not=j_0}D_j$ is located in the half-space $x\cdot\omega<h_D(\omega)-\delta$.
Therefore, we have
$$\displaystyle
e^{-\tau h_D(\omega)}\sum_j(k_j-1)
\int_{\partial D_j}u\frac{\partial v}{\partial\nu}
=e^{-\tau h_{D_{j_0}}(\omega)}(k_{j_0}-1)
\int_{\partial D_{j_0}}u\frac{\partial v}{\partial\nu}+O(\tau e^{-\tau\delta})
$$
as $\tau\longrightarrow\infty$.

Using a well known expansion of $u$ about $x_0$ (see, for instance, Proposition 2.1 in \cite{Ie2})
and (3.6) for $c=u(x_0)$, we obtain the asymptotic exapnsion as $\tau\longrightarrow\infty$:
$$\displaystyle
e^{-\tau h_{D_{j_0}}(\omega)}(k_{j_0}-1)\int_{\partial D_{j_0}}u\frac{\partial v}{\partial\nu}
\sim
e^{i\tau x_0\cdot\omega^{\perp}}\sum_{j=1}^{\infty}\frac{L_j}{\tau^{\mu_j}}
$$
where $0<\mu_1<\mu_2<\cdots$ (Proposition 3.2 in \cite{Ie2}).  Then the problem is to show that
$L_j\not=0$ for some $j$.
We see that if $L_j=0$ for all $j$, then $u$ has a harmonic continuation in a neighbourhood of $x_0$
having a rotation invariance property with respect to some angle $0<\theta\le\pi$ (Lemma 4.1 in \cite{Ie2}).
Note that $x_0$ is also a vertex of the convex hull of $D$ which is a polygon.
Then applying Friedman-Isakov's extension argument \cite{FI} to $u$ outside the convex hull of $D$,
we see that $u$ has a harmonic continuation in the whole domain.
Then it is easy to see that $u$ has to be a constant: a contradiction.

\noindent
$\Box$

From (2.2) we see that $\mbox{$\cal D$}=\mbox{$\cal D$}(P,Q;\,\cdot\,)\in H^1_{\text{loc}}(\Omega)$ satisfies
$$\displaystyle
\int_{\Omega}\gamma\nabla\mbox{$\cal D$}\cdot\nabla\varphi dx=0
$$
for all $\varphi\in C^{\infty}_0(\Omega)$; $\mbox{$\cal D$}(P,Q;\,\cdot\,)$ is not a constant function.
Therefore, from Lemma 3.1 for $u=\mbox{$\cal D$}(P,Q;\,\cdot\,)$, (1.5) and (3.1) we obtain the two formulae in Theorem 1.1.

$$\quad$$

\centerline{{\bf Acknowledgments}}

The author thanks the referees for several suggestions for the improvement
of the manuscript.
This research was partially supported by Grant-in-Aid for
Scientific Research (C)(No. 11640151) of Japan  Society for
the Promotion of Science.

$$\quad$$


\begin{thebibliography}{99}



\bibitem{AD}  Alessandrini, G. and DiBenedetto, E.,
              Determining $2$-dimensional cracks in $3$-dimensional bodies: uniqueness and stability,
              Indiana Univ. Math. J., {\bf 46}(1997), 1-82.




\bibitem{C} Calder\'on,  A. P.,
       \newblock On an inverse boundary value problem,
       \newblock in Seminar on Numerical Analysis and its Applications to Continuum Physics
       (Meyer, W. H. and Raupp, M. A. eds.), Brazilian Math. Society, Rio de Janeiro, 1980,
        65-73.




\bibitem{FI} Friedman, A. and Isakov, M.,
         \newblock On the uniqueness in the inverse conductivity problem with one measurements,
         \newblock Indiana Univ. Math. J., {\bf 38}(1989), 563-579.




\bibitem{Gr} Grisvard, P.,
          \newblock  Elliptic problems in nonsmooth domains, Pitman, Boston, 1985.




\bibitem{Ie}  Ikehata, M.,
              Enclosing a polygonal cavity in a two-dimensional bounded domain from Cauchy data,
              Inverse Problems, {\bf 15}(1999), 1231-1241.



\bibitem{Ie2}  Ikehata, M.,
               On reconstruction in the inverse conductivity problem with one measurement,
               Inverse Problems, {\bf 16}(2000), 785-793.





\bibitem{IS} Isakov, V.,
            \newblock On uniqueness of recovery of a discontinuous conductivity coefficients,
            \newblock Comm. Pure. Appl. Math., {\bf 41}(1988), 865-877.




              


\bibitem{L}  Ladyzhenskaya, O. A. and Ural'tzeva N. N.,
             \newblock Linear and quasilinear elliptic equations, 1968, London, Academic Press.







\bibitem{N2} Nachman, A.,
             Global uniqueness for a two-dimensional inverse boundary value problem,
             Ann. of Math., {\bf 143}(1996), 71-96.






\bibitem{U} Uhlmann, G.,
    \newblock  Developments in inverse problems since Calder\'on's foundational paper,
    \newblock in Harmonic analysis and partial differential equations
    (Christ, M., Kenig, C. E. and Sadosky, C., eds.), 1999, pp. 295-345,
    The University of Chicago Press, Chicago and London.






\end{thebibliography}
\end{document}